\documentclass[12pt,reqno]{amsart}
\usepackage[utf8]{inputenc}
\usepackage{amsmath, amsfonts, amssymb, amsthm, hyperref}
\usepackage{bm}
\allowdisplaybreaks[4]
\def\l{\left}
\def\r{\right}
\def\bg{\bigg}
\def\({\bg(}
\def\){\bg)}
\def\t{\text}
\def\f{\frac}

\def\sm{\setminus}

\def\eq{\equiv}

\def\1{{\bf 1}}

\def\pmod #1{\ ({\rm{mod}}\ #1)}

\def\<{\langle}
\def\>{\rangle}
\def\t{\text}

\theoremstyle{plain}
\newtheorem{theorem}{Theorem}[section]
\newtheorem{lemma}{Lemma}
\newtheorem{corollary}{Corollary}
\newtheorem{conjecture}{Conjecture}

\theoremstyle{definition}
\newtheorem*{Ack}{Acknowledgment}

\theoremstyle{remark}

\makeatletter
\@namedef{subjclassname@2020}{%
  \textup{2020} Mathematics Subject Classification}
\makeatother
 \vspace{4mm}

\begin{document}

\medskip

\title[On certain determinants and the square roots of some determinants involving  Legendre Symbols]
{On certain determinants and the square roots
\\of some determinants involving  Legendre Symbols }
\author{Chen-kai Ren}
\address {(Chen-kai Ren) Department of Mathematics, Nanjing
University, Nanjing 210093, People's Republic of China}
\email{ckren@smail.nju.edu.cn}
\author{Xin-qi Luo}
\address {(Xin-qi Luo ) Department of Mathematics, Nanjing
University, Nanjing 210093, People's Republic of China}
\email{lxq15995710087@163.com}

\keywords{determinant, Legendre Symbol.
\newline \indent 2020 {\it Mathematics Subject Classification}. Primary 11A25; Secondary 11N25.
\newline \indent Supported by the Natural Science Foundation of China (grant no. 12371004).}
\begin{abstract}
Let $p>3$ be a prime and $(\frac{.}{p})$ be the Legendre symbol. For any integer $d$ with $p\nmid d$ and any positive integer $m$, Sun introduced the determinants
$$T_m(d,p)=\det\left[(i^2+dj^2)^m\l(\frac{i^2+dj^2}{p}\r)\right]_{1\leqslant i,j \leqslant (p-1)/2},$$
and
$$ D_p^{(m)}= \det\l[(i^2-j^2)^m\l(\frac{i^2-j^2}{p}\r)\r]_{1\leq i,j\leq (p-1)/2} .$$
In this paper, we obtain some properties of $T_m
(d,p)$ and $ \sqrt{D_p^{(m)}}$ for some $m$.   We also  confirm
some related conjectures posed by Zhi-Wei Sun.

\end{abstract}
\maketitle

\section{Introduction}
\setcounter{lemma}{0}
\setcounter{theorem}{0}
\setcounter{equation}{0}
\setcounter{conjecture}{0}
\setcounter{remark}{0}
\setcounter{corollary}{0}
Determinants of certain special matrices are useful in many branches of mathematics. Readers may refer to Krattenthaler’s survey paper \cite{K04} for recent progress and advanced techniques on this topic. In this paper, we study some determinants involving Legendre symbols.

Let $p$ be an odd prime, and let $(\frac{.}{p})$ be the Legendre symbol.
Carlitz \cite{C} determined the characteristic polynomial of the matrix
$$\left[x+\left(\frac{i-j}{p}\right)\right]_{1\leqslant i,j\leqslant p-1},
$$
and Chapman \cite{C1} used quadratic Gauss sums to determine the values of the determinant
$$
\det\left[ x+\left(\frac{i+j-1}{p}\right)\right]_{1\leqslant i,j\leqslant (p-1)/{2}}.
$$
Vsemirnov \cite{V12,V13} confirmed a challenging conjecture of Chapman
by evaluating the determinant
$$\det\left[\left(\frac{j-i}{p}\right)\right]_{1\leqslant i,j\leqslant (p+1)/{2}}.
$$
Let $d$ be any integer. Sun \cite{S19} introduced the determinant
$$ S(d,p)=\det\l[\l(\frac{i^2+dj^2}{p}\r)\r]_{1\leqslant i,j \leqslant (p-1)/2}.
$$
Sun \cite[Theorem1.2]{S19} proved that
$$
   \l(\frac{S(d,p)}{p}\r)=
\begin{cases}
   (\frac{-1}{p})  & \t{if} \ (\frac{d}{p})=1, \\
    \quad 0  & \t{if} \ (\frac{d}{p})=-1.\\
\end{cases}
$$
D.Grinberg, Sun and L.Zhao \cite{GSZ} showed that if $p>3$
then
$$ \det\l[(i^2+dj^2)\l(\frac{i^2+dj^2}{p}\r)\r]_{0\leqslant i,j \leqslant (p-1)/2}\equiv 0 \pmod p.
$$
Motivated by this, for any positive integer $n$ with $ (p-1)/2 \leq n \leq p-1$ Sun introduced the determinant
$$S_n(d,p)=\det\left[(i^2+dj^2)^{n}\right]_{1\leqslant i,j \leqslant (p-1)/2}.
$$
In 2022, H.-L. Wu, Y.-F. She and L.-Y. Wang \cite{WSW} proved the conjecture of Sun \cite[Conjecture 4.5]{S19} that if $p>3$ and $p \nmid d$ then
$$
\l(\frac{S_{(p+1)/2}(d,p)}{p}\r)=
\begin{cases}
     (\frac{d}{p})^{(p-1)/4}  & \t{if} \ p\equiv1 \pmod 4, \\
    (\frac{d}{p})^{(p+1)/4}(-1)^{(h(-p)-1)/2}  & \t{if}\ p\equiv 3 \pmod 4,\\
\end{cases}
$$
where $h(-p)$ denotes the class number of the imaginary quadratic field $\mathbb Q(\sqrt{-p})$.

For any prime $p\eq3\pmod4$, Sun \cite[Remark 1.3]{S19} showed that
$$ S_{p-2}(1,p)= \det \left[\frac{1}{i^2+j^2}\right]_{1\leqslant i,j \leqslant (p-1)/2}\eq\l(\f 2p\r)\pmod p.
$$
Sun \cite{S24} proved that if $p\equiv 3\pmod 4, $

 $$\l(\frac{S_{p-3}(1,p)}{p}\r)=1.$$
 Ren and Sun \cite{R24} solved some conjectures of Sun concerning the Legendre symbols of $S_{p-2},$ $S_{p-3}$ and $S_{p-4}$.
Sun posed the following conjecture \cite[Remark 1.1]{S24}.
 \begin{conjecture} Let $p>3$ be a prime with $p\neq 11$. Let
 $$T_2(1,p)=\det \l[{(i^2+j^2)^2\l(\frac {i^2+j^2}p\r)}\r]_{1\leq i,j\leq (p-1)/2},$$
where $(\frac{.}{p})$ is the Legendre symbol. Then
$$\l(\f{T_2(1,p)}p\r)=\begin{cases} \l(\frac {2} p\r)& \t{if}\ p\equiv 1\pmod 4,
\\\l(\frac {-6}{p}\r)& \t{if}\ p\equiv 3\pmod 4.\end{cases}$$
\end{conjecture}
In this paper, we confirm this conjecture. Inspired by this, we start to study determinants on finite field.
Let $q=2n+1$ be an odd prime power and let $\phi$ denote the unique quadratic multiplicative character of $\mathbb F_q$,
which is the map $\mathbb{F}_q \to \mathbb{C} $ that sends $0$ to $ 0$, each non-zero square to $1$, and each non-square to $-1$. Let $\mathbb F_q^{\times}:=\mathbb F_q\sm \{0\}$ and let
$$
\mathbb F_q^{\times 2}:=\{x^2:x\in \mathbb F_q^{\times}\}=\{a_1, a_2,\cdots a_n\}.
$$
We define
$$\widetilde T_m(d,q)=\det[ {(a_i+da_j)^m\phi(a_i+da_j)}]_{1\leq i,j\leq n}.$$
We have  the following generalized result.
\begin{theorem}
Let $q=2n+1$ be an odd prime power and let $d\in \mathbb F_q^{\times} $. Suppose $\mathrm{char}(\mathbb F_q)>3$. Let
$$\widetilde T_2(d,q)=\det[ {(a_i+da_j)^2\phi(a_i+da_j)}]_{1\leq i,j\leq n}.$$
If we view $\widetilde T_2(d,q)$ as a determinant over $\mathbb F_q$, then the following results hold:

(\romannumeral1) If $q\equiv 1\pmod 4$, then
$$\widetilde T_2(d,q)=d^{\frac {q-1}{4}}\times \frac {q-1}2!x_q(d)^2$$
for some $x_q(d)\in \mathbb F_q$.

(\romannumeral2) If $q\equiv 3\pmod 4$, then
$$\widetilde T_2(d,q)=(-1)^{\frac {q-3}4}3d^{\frac {q+1}{4}}y_q(d)^2$$
for some $y_q(d)\in \mathbb F_q$.
\end{theorem}
Applying this result, we can obtain Conjecture 1.1.
\begin{corollary}
    Conjecture 1.1 holds.
\end{corollary}

\begin{theorem}
Let $q=2n+1$ be an odd prime power and let $d\in \mathbb F_q^{\times}$. Suppose $\mathrm{char}(\mathbb F_q)>7$. Let
$$\widetilde T_{(q-11)/2}(d,q)=\det[ {(a_i+da_j)^{(q-11)/2}\phi(a_i+da_j)}]_{1\leq i,j\leq n}.$$
If we view $\widetilde T_{(q-11)/2}(d,q)$ as a determinant over $\mathbb F_q$ , then the following results hold:

(\romannumeral1) If $q\equiv 1\pmod 4$, then
$$\widetilde T_{(q-11)/2}(d,q)=d^{\frac{q-1}{4}}x_q(d)^2$$
   for some $x_q(d)\in \mathbb F_q$.

(\romannumeral2) If $q\equiv 3\pmod 4$, then
$$\widetilde T_{(q-11)/2}(d,q)=7(q-1)!d^{\frac{q+1}{4}}y_q(d)^2$$
for some $y_q(d)\in \mathbb F_q$.
\end{theorem}
Applying this result, we can get the following corollary.
\begin{corollary}Let $p>7$ be a prime and let $d$ be any integer with $\l(\frac{d}{p}\r)=1$. Let
$$ T_{(p-11)/2}(d,p)=\det \l[(i^2+dj^2)^{(p-11)/2}\l(\frac{i^2+dj^2}{p}\r)\r]_{1\leq i,j\leq  (p-1)/2}.$$
If $\(\frac{T_{(p-11)/2}(d,p)}{p}\)= -1 $, then $p \equiv  1\ \t{or}\   2\   \t{or}\  4 \pmod 7$. \end{corollary}
Let $m,n\in \mathbb Z^{+}$ with $n$ odd. Sun introduced the determinant
$$ D_n^{(m)}= \det\l[(i^2-j^2)^m\l(\frac{i^2-j^2}{n}\r)\r]_{1\leq i,j\leq (n-1)/2} .$$
If $2\nmid m$ and $4\mid n-1$, then $D_n^{(m)}$
is skew-symmetric and of even order. Hence, it is the square of a certain integere by
Cayley’s Theorem \cite{Sta}.
We also prove the following conjectures of Sun \cite{S24}.
\begin{theorem}
For any positive odd integer $m$, the set
$$
E(m)=\{p: p\ is \ a \ prime\ with\ 4\mid p-1\ and\ p\mid D_p^{(m)}\}
$$
is finite. In particular,
$$E(5)=\{29\},\ E(7)=\{13, 53\},\ E(9)=\{13, 17, 29\} $$
$$E(11)=\{17, 29\} \ \t{and} \ E(13)=\{17, 109, 401\}. $$
\end{theorem}
\begin{theorem}
For any prime $p\equiv 1\pmod 4,$ we have
$$
\l(\frac{\sqrt{D_p^{(1)}}}{p}\r)=(-1)^{|\{0<k<\frac p4: (\frac kp)=-1\}|}\l(\frac p3\r).
$$
\end{theorem}
\begin{theorem}
For any prime $p\equiv 1\pmod 4,$ we have
$$
\l(\frac{\sqrt{D_p^{(3)}}}{p}\r)=(-1)^{|\{0<k<\frac p4: (\frac kp)=-1\}|}\l(\frac {p}{4+(-1)^{(p-1)/4}}\r).
$$
\end{theorem}
\section{Some  Preparations before the proof}
\setcounter{lemma}{0}
\setcounter{theorem}{0}
\setcounter{equation}{0}
\setcounter{conjecture}{0}
\setcounter{remark}{0}
\setcounter{corollary}{0}

Let $q=2n+1$ be an ood prime power and  recall $\phi$ denote the unique quadratic multiplicative characteristic $\mathbb F_q$. Let $\mathbb F_q^{\times}:=\mathbb F_q\sm \{0\}$ and let
$$
\mathbb F_q^{\times 2}:=\{x^2:x\in \mathbb F_q^{\times}\}=\{a_1, a_2,\cdots a_n\}.
$$
For any $x,y\in \mathbb{F}_q$, we write $x\equiv y \pmod{\mathbb{F}_q^{\times 2}}$ if there is an element $z\in \mathbb{F}_q^{\times }$  such that $x=yz^2.$

We need the following  lemma \cite{K04} on determinants.

\begin{lemma}\label{lemma4}
Let $R$ be a commutative ring with identity, and let $P(x)=\sum_{i=0}^{n-1}b_ix^i\in R[x]$. Then we have
$$\det{\left[P(X_iY_j)\right]}_{1\leqslant i<j\leqslant n}=b_0b_1\cdots b_{n-1} \prod_{1\leqslant i< j\leqslant n}(X_i-X_j)(Y_i-Y_j).$$
\end{lemma}
Let  $a$ be an integer with $(a,n)=1.$ We notice that multiplication by $a$ introduces a permutation $\tau_a$ of $\mathbb Z/n\mathbb Z$. Lerch \cite{L96} obtained the following result which determines the sign of $\tau_a$.
\begin{lemma}
Let $\mathrm{sgn}(\tau_a^{(n)})$ denote the sign of $\tau_a$. Then
$$\mathrm{sgn}(\tau_a^{(n)})=\begin{cases} (\frac an)&\t{if}\ n\equiv 1\pmod 2,
\\1&\t{if}\ n\equiv 2\pmod 4,
\\(-1)^{\frac {a-1}2}&\t{if}\ n\equiv 0\pmod 4,\end{cases}$$
where $(\frac {\cdot }n)$ denotes Jacobi symbol if $n$ is odd.
\end{lemma}
Clearly $\mathrm{inv_q}:a_i\to a_i^{-1}$ is a permutation of $a_1,\cdots, a_n$. Fix a generator $g$ of the cyclic group $\mathbb F_q^{\times }$. If we reindex the elements $g^0,g^1,\cdots ,g^{2(n-1)}$ of $\mathbb F_q^{\times^2}$ as $0,1,\cdots, n-1$ of $\mathbb Z/n\mathbb Z$, then the permutation $\mathrm{inv_q}$ becomes the permutation $\tau_{-1}$ that sends $i$ to $n-i$ for $1\leq i\leq n-1$ while leaving $0$ fixed. This permutation has $\frac {(n-1)(n-2)}2$ inversions. In view of the above we obtain the following result.
\begin{lemma}
Let notations be as above. Then
$$
\mathrm{sgn}(\mathrm{inv_q})=(-1)^{\frac {(n-1)(n-2)}2}=(-1)^{\frac {(q-3)(q-5)}8}.
$$
\end{lemma}
\section{proofs of Theorem 1.1 and 1.2}
\setcounter{lemma}{0}
\setcounter{theorem}{0}
\setcounter{equation}{0}
\setcounter{conjecture}{0}
\setcounter{remark}{0}
\setcounter{corollary}{0}

\noindent{\bf Proof of Theorem 1.1}.
Recall that $\phi $ denote the unique quadratic multiplicative character of $F_q$.
If we view $\widetilde T_2(d,q)$ as a determinant over $\mathbb F_q$, then we have

$$\widetilde T_2(d,q)=\det \left[\l(a_i+da_j\r)^{\frac{q+3}{2}}\right]_{1\leq i,j\leq n}.$$
Thus
$$\widetilde T_2(d,q)\equiv \det \left[\l(\frac {a_i}{a_j}+d\r)^{\frac{q+3}{2}}\right]_{1\leq i,j\leq n}\pmod {\mathbb F_q^{\times 2}}.$$
Noting that
$$
\l(\frac {a_i}{a_j}\r)^{\frac {q-1}2}=1 \ \t{and} \ d^{\frac {q-1}2}=\phi (d),
$$
we can verify that

\begin{align*}
&\l(\frac {a_i}{a_j}+d\r)^{\frac{q+3}{2}}\\
=&\sum_{k=0}^{\frac{q+3}{2}}\binom{\frac{q+3}{2}}{k}\l(\frac {a_i}{a_j}\r)^kd^{\frac{q+3}{2}-k}\\
=&\sum_{k=3}^{\frac{q-3}{2}}\binom{\frac{q+3}{2}}{k}\l(\frac {a_i}{a_j}\r)^kd^{\frac{q+3}{2}-k}+d^{\frac{q+3}{2}}+\frac{q+3}{2}\l(\frac {a_i}{a_j}\r)d^{\frac{q+1}{2}}+\frac{(q+3)(q+1)}{8}d^{\frac{q-1}{2}}\l(\frac {a_i}{a_j}\r)^2\\
&+\frac{(q+3)(q+1)}{8}d^{2}\l(\frac {a_i}{a_j}\r)^{\frac {q-1}2}+\frac{q+3}{2}d\l(\frac {a_i}{a_j}\r)^{\frac {q+1}2}+\l(\frac {a_i}{a_j}\r)^{\frac {q+3}2}\\
=&\sum_{k=3}^{\frac{q-3}{2}}\binom{\frac{q+3}{2}}{k}\l(\frac {a_i}{a_j}\r)^kd^{\frac{q+3}{2}-k}+\l(\frac{3}{8}+\phi (d)\r)d^2+\frac{3}{2}\l(\frac {a_i}{a_j}\r)d(\phi (d)+1)+\l(\frac {a_i}{a_j}\r)^2\l(\frac{3}{8}\phi (d)+1\r).\\
\end{align*}
Hence
$$
\l(\frac {a_i}{a_j}+d \r)^{\frac{q+3}{2}}=f\l(\frac {a_i}{a_j}\r)
$$
where

\begin{align*}
f(T)=&\l(\frac{3}{8}+\phi (d)\r)d^2+\frac{3}{2}d(\phi (d)+1)T+\l(\frac{3}{8}\phi (d)+1\r)T^2\\
&+\sum_{k=3}^{\frac{q-3}{2}}\binom{\frac{q+3}{2}}{k}T^kd^{\frac{q+3}{2}-k}.
\end{align*}
Observe that
$$
\l(\frac{3}{8}+\phi (d)\r)\l(\frac{3}{8}\phi (d)+1\r)\equiv 0 \ \t{or} \ \phi(d)\pmod {\mathbb F_q^{\times^2}}
$$
(just check both cases $\phi (d)=1$ or $\phi (d)=-1$ and consider whether $\frac{3}{8}+\phi (d)\equiv 0 \pmod q$  ).
Let $C_f$ denote the product of coefficients of $f(T),$
then
\begin{align*}
C_f=&(\frac{3}{8}+\phi (d))\frac{3}{2}d^3(\phi (d)+1)(\frac{3}{8}\phi (d)+1)\prod_{3}^{\frac{q-3}{2}}\binom{\frac{q+3}{2}}{k}d^{\frac{q+3}{2}-k}\\
=&\frac{3}{2}d^3(\phi (d)+1)(\frac{3}{8}+\phi (d))(\frac{3}{8}\phi (d)+1)\prod_{3}^{\frac{q-3}{2}}\binom{\frac{q+3}{2}}{k}\times d^{\frac {(q+3)(q-7)}{8}}\\
\equiv &6d^{\frac {(q-3)(q-1)}{8}}(\phi (d)+1)(\frac{3}{8}+\phi (d))(\frac{3}{8}\phi (d)+1)\prod_{3}^{\frac{q-3}{2}}\binom{\frac{q+3}{2}}{k}\pmod {\mathbb F_q^{\times^2}}.
\end{align*}
Thus applying Lemma 2.1 to $P(T)=f(T),$ we rewrite as
\begin{align*}
\widetilde T_2(d,q)&\equiv \mathrm{sgn}(\mathrm{inv_q})6d^{\frac {(q-3)(q-1)}{8}}(\phi (d)+1)(\frac{3}{8}+\phi (d))(\frac{3}{8}\phi (d)+1)\prod_{3}^{\frac{q-3}{2}}\binom{\frac{q+3}{2}}{k}\prod_{1\leq i<j\leq n}(a_j-a_i)^2\\
&\equiv \mathrm{sgn}(\mathrm{inv_q})6d^{\frac {(q-3)(q-1)}{8}}(\phi (d)+1)(\frac{3}{8}+\phi (d))(\frac{3}{8}\phi (d)+1)\prod_{3}^{\frac{q-3}{2}}\binom{\frac{q+3}{2}}{k} \pmod {\mathbb F_q^{\times^2}} .\\
\end{align*}
Now we divide the remaining proof into two cases.

 \textbf{Case 1.} $q\equiv 1\pmod 4$.

In this case, since $\pm 1\in {\mathbb F_q^{\times 2}}$ and $\mathrm{char}(\mathbb F_q)> 3$, we can verify that
$$
\mathrm{sgn}(\mathrm{inv_q})6d^{\frac {(q-3)(q-1)}{8}}\equiv 6d^{\frac {q-1}{4}}\pmod {\mathbb F_q^{\times 2}}
$$
and that
$$
\prod_{k=3}^{\frac{q-3}{2}}\binom{\frac{q+3}{2}}{k}\in \{0\}\cup \binom{\frac{q+3}{2}}{\frac{q+3}{4}}\mathbb F_q^{\times^2}
$$
where $\binom{\frac{q+3}{2}}{\frac{q+3}{4}}\equiv 3\times \frac {q-1}2!\pmod {\mathbb F_q^{\times 2}}.$
Hence, we have
\begin{align*}
 \widetilde T_2(d,q)&\equiv 6d^{\frac {q-1}{4}}(\phi (d)+1)\phi(d)\times 3\frac {q-1}2! \ or \ 0  \\
 &\equiv d^{\frac {q-1}{4}}\times \frac {q-1}2!  \ \t{or} \ 0  \pmod {\mathbb F_q^{\times^2}}\\
\end{align*}
Therefore, there is an element $x_q(d)\in \mathbb F_q$ such that
$$
\widetilde T_2(d,q)=d^{\frac {q-1}{4}}\times \frac {q-1}2!x_q(d)^2.
$$

 \textbf{Case 2. }$q\equiv 3\pmod 4.$

In this case, we have
$$\prod_{k=3}^{\frac{q-3}{2}}\binom{\frac{q+3}{2}}{k}\in \{0\}\cup \mathbb F_q^{\times^2}.$$
By Lemma 2.3, we obtain that
$$
\mathrm{sgn}(\mathrm{inv_q})6d^{\frac {(q-3)(q-1)}{8}}\equiv (-1)^{\frac {q-3}4}6d^{\frac {q-3}{4}}\pmod {\mathbb F_q^{\times^2}}.
$$
Noting that $\phi(d)\equiv d  \pmod {\mathbb F_q^{\times^2}},$ we verify that
\begin{align*}
 \widetilde T_2(d,q)&\equiv  (-1)^{\frac {q-3}4}6d^{\frac {q-3}{4}}(\phi (d)+1)\phi(d)  \ \t{or} \ 0  \\
 & \equiv(-1)^{\frac {q-3}4}3d^{\frac {q+1}{4}}  \ \t{or} \ 0  \pmod {\mathbb F_q^{\times^2}}.\\
\end{align*}
In view of the above there is an element $y_q(d)\in \mathbb F_q$ such that
$$
\widetilde T_2(d,q)=(-1)^{\frac {q-3}4}3d^{\frac {q+1}{4}}y_q(d)^2.
$$  \qed

\noindent{\bf Proof of Corollary 1.1}.
If we view $T_2(1,p)$ as the determinant  over $\mathbb F_p$, we can take advantage of the proof above. Take $d=1$ and let $q=p$ with  $p\geq 5$ and $p\neq 11$. We obtain that

\begin{align*}
T_2(1,p)&\equiv 3\times\mathrm{sgn}(\mathrm{inv_p})(\frac{3}{8}+1 )^2\prod_{3}^{\frac{p-3}{2}}\binom{\frac{p+3}{2}}{k} \\
&\equiv 3\times \mathrm{sgn}(\mathrm{inv_p})\prod_{3}^{\frac{p-3}{2}}\binom{\frac{p+3}{2}}{k} \pmod {\mathbb F_p^{\times^2}}.\\
\end{align*}
Also, we divide the remaining proof into two cases.

  \textbf{Case 1.} $p\equiv 1\pmod 4.$ In this case, we have
 \begin{align*}
     T_2(1,p) &\equiv 3\binom{\frac{p+3}{2}}{\frac{p+3}{4}}\prod_{3}^{\frac{p-1}{4}}\binom{\frac{p+3}{2}}{k}^2 \\
     &\equiv \frac {p-1}2! \pmod {\mathbb F_p^{\times^2}}.\\
 \end{align*}
Hence, by  \cite[Lemma 2.3]{S19} , we have
$$\(\frac{T_2(1,p)}{p}\)= \(\frac {(\frac {p-1}2)!} p\)=\(\frac {2} p\).$$

\textbf{ Case 2.} $p\equiv 3\pmod 4.$ In this case, by Lemma 2.3 we have
 \begin{align*}
     T_2(1,p) &\equiv 3\times\mathrm{sgn}(\mathrm{inv}_q)\prod_{3}^{\frac{p+1}{4}}\binom{\frac{p+3}{2}}{k}^2 \\
     &\equiv 3(-1)^{\frac {p-3}4} \pmod {\mathbb F_p^{\times^2}}.\\
 \end{align*}
Therefore, for $\l(\frac {(-1)^{\frac {p-3}4}}p\r)=(-1)^{\frac {p-3}4}=-1\times (-1)^{\frac {(p-1)(p+1)}8}=\l(\frac {-2}p\r), $ we obtain
\begin{align*}
\l(\frac {T_2(1,p)}p\r) &= \l(\frac{3(-1)^{\frac {p-3}4}}{p}\r)=\l(\frac {-6}p\r).\\
\end{align*}  \qed

\noindent{\bf Proof of Theorem 1.2}. Recall that $\phi $ denote the unique quadratic multiplicative character of $\mathbb F_q$.
If we view $\widetilde T_{(q-11)/2}(d,q)$ as a determinant over $\mathbb F_q$, then we have

$$\widetilde T_{(q-11)/2}(d,q)=\det \left[(a_i+da_j)^{q-6}\right]_{1\leq i,j\leq n}.$$
Thus
$$\widetilde T_{(q-11)/2}(d,q)\equiv \det \left[\l(\frac {a_i}{a_j}+d\r)^{q-6}\right]_{1\leq i,j\leq n}\pmod {\mathbb F_q^{\times^2}}.$$
Noting that
$$
\l(\frac {a_i}{a_j}\r)^{\frac {q-1}2}=1 \ and \ d^{\frac {q-1}2}=\phi (d),
$$
we can verify that
\begin{align*}
&\l(\frac {a_i}{a_j}+d\r)^{q-6}\\
=&\sum_{k=0}^{q-6}\binom{q-6}{k}\l(\frac {a_i}{a_j}\r)^kd^{q-6-k}\\
=&\sum_{k=0}^{\frac{q-11}{2}}\l(\binom{q-6}{k}d^{q-6-k}+\binom{q-6}{\frac{q-1}{2}+k}d^{q-6-k-(q-1)/2}\r)\l(\frac {a_i}{a_j}\r)^k \\
&+\sum_{k=0}^{3}\binom{q-6}{\frac{q-9}{2}+k}d^{(q-3)/2-k}\l(\frac {a_i}{a_j}\r)^{(q-9)/2+k}\\
=&\sum_{k=0}^{\frac{q-11}{2}}\l(\binom{q-6}{k}+\phi(d)^{-1}\binom{q-6}{\frac{q-11}{2}-k}\r)d^{q-6-k}\l(\frac {a_i}{a_j}\r)^k \\
&+\sum_{k=0}^{3}\binom{q-6}{\frac{q-9}{2}+k}d^{(q-3)/2-k}\l(\frac {a_i}{a_j}\r)^{(q-9)/2+k}.
\end{align*}
Hence
$$
\l(\frac {a_i}{a_j}+d \r)^{q-6}=f\l(\frac {a_i}{a_j}\r)
$$
where

\begin{align*}
f(T)=&\sum_{k=0}^{\frac{q-11}{2}}\l(\binom{q-6}{k}+\phi(d)^{-1}\binom{q-6}{\frac{q-11}{2}-k}\r)d^{q-6-k}T^k \\
&+\sum_{k=0}^{3}\binom{q-6}{\frac{q-9}{2}+k}d^{(q-3)/2-k}T^{(q-9)/2+k}.
\end{align*}
Observe that
$$
\l(\binom{q-6}{k}+\phi(d)^{-1}\binom{q-6}{\frac{q-11}{2}-k}\r)\l(\binom{q-6}{\frac{q-11}{2}-k}+\phi(d)^{-1}\binom{q-6}{k}\r)\equiv 0\ or\ \phi(d)\pmod {\mathbb F_q^{\times^2}}
$$
(just check both cases $\phi (d)=1$ or $\phi (d)=-1$ and consider whether $\binom{q-6}{k}+\phi(d)^{-1}\binom{q-6}{\frac{q-11}{2}-k}\equiv 0 \pmod q$ ).
Let $C_f$ denote the product of coefficients of $f(T),$
then
\begin{align*}
C_f=&\prod_{k=0}^{3}\binom{q-6}{\frac{q-9}{2}+k}\prod_{k=0}^{3}d^{(q-3)/2-k}
\prod_{k=0}^{\frac{q-11}{2}}\l(\binom{q-6}{k}+\phi(d)^{-1}\binom{q-6}{\frac{q-11}{2}-k}\r)\prod_{k=0}^{\frac{q-11}{2}}d^{q-6-k}\\
= &\binom{q-6}{\frac{q-9}{2}}^2\binom{q-6}{\frac{q-7}{2}}^2 \prod_{k=0}^{\frac{q-11}{2}}\l(\binom{q-6}{k}+\phi(d)^{-1}\binom{q-6}{\frac{q-11}{2}-k}\r)d^{\frac{(q-1)(3q-21)}{8}}.\\
\end{align*}
Applying Lemma 2.1 to $P(T)=f(T),$ we rewrite as
\begin{align*}
\widetilde T_{(q-11)/2}(d,q)
\equiv& \mathrm{sgn}(\mathrm{inv_q})\prod_{k=0}^{\frac{q-11}{2}}\l(\binom{q-6}{k}+\phi(d)^{-1}\binom{q-6}{\frac{q-11}{2}-k}\r)d^{\frac{(q-1)(3q-21)}{8}} \\
&\times\binom{q-6}{\frac{q-9}{2}}^2\binom{q-6}{\frac{q-7}{2}}^2 \prod_{1\leq i<j\leq n}(a_j-a_i)^2\\
\equiv& \mathrm{sgn}(\mathrm{inv_q})\prod_{k=0}^{\frac{q-11}{2}}\l(\binom{q-6}{k}+\phi(d)^{-1}\binom{q-6}{\frac{q-11}{2}-k}\r) \\
&\times\binom{q-6}{\frac{q-9}{2}}^2\binom{q-6}{\frac{q-7}{2}}^2d^{\frac{(q-1)(3q-21)}{8}}\pmod {\mathbb F_q^{\times^2}}.\\
\end{align*}
Now we divide the remaining proof into two cases.

 \textbf{Case 1.} $q\equiv 1\pmod 4$.

Since $\pm 1\in \mathbb{F}_q^{\times^2}$ and $\mathrm{char}(\mathbb{F}_q)>7$, we can verify that
\begin{align*}
\widetilde T_{(q-11)/2}(d,q) &\equiv\mathrm{sgn}(\mathrm{inv_q})\prod_{k=0}^{\frac{q-11}{2}}\l(\binom{q-6}{k}+\phi(d)^{-1}\binom{q-6}{\frac{q-11}{2}-k}\r)\binom{q-6}{\frac{q-9}{2}}^2\binom{q-6}{\frac{q-7}{2}}^2d^{\frac{(q-1)(3q-21)}{8}} \\
 &\equiv \prod_{k=0}^{\frac{q-13}{4}}\l(\binom{q-6}{k}+\phi(d)^{-1}\binom{q-6}{\frac{q-11}{2}-k}\r)\l(\binom{q-6}{\frac{q-11}{2}-k}+\phi(d)^{-1}\binom{q-6}{k}\r)\\
 &\times\binom{q-6}{\frac{q-9}{2}}^2\binom{q-6}{\frac{q-7}{2}}^2d^{\frac{(q-1)(3q-21)}{8}} \\
 &\equiv d^{\frac{(q-1)(3q-21)}{8}}\phi(d)^{\frac{q-9}{4}} \ or \ 0  \\
&\equiv d^{\frac{q-1}{4}} \ or \ 0 \pmod {\mathbb F_q^{\times^2}}. \\
\end{align*}
Hence there is an element $x_q(d)\in \mathbb{F}_q$ such that
$$
\widetilde T_{(q-11)/2}(d,q)=d^{\frac{q-1}{4}}x_q(d)^2.
$$

 \textbf{Case 2. }$q\equiv 3\pmod 4.$
We can verify that
\begin{align*}
  \widetilde T_{(q-11)/2}(d,q)&\equiv\mathrm{sgn}(\mathrm{inv_q})\prod_{k=0}^{\frac{q-11}{2}}\l(\binom{q-6}{k}+\phi(d)^{-1}\binom{q-6}{\frac{q-11}{2}-k}\r)\binom{q-6}{\frac{q-9}{2}}^2\binom{q-6}{\frac{q-7}{2}}^2d^{\frac{(q-1)(3q-21)}{8}} \\
 &\equiv \mathrm{sgn}(\mathrm{inv_q})\prod_{k=0}^{\frac{q-15}{4}}\l(\binom{q-6}{k}+\phi(d)^{-1}\binom{q-6}{\frac{q-11}{2}-k}\r)\l(\binom{q-6}{\frac{q-11}{2}-k}+\phi(d)^{-1}\binom{q-6}{k}\r)\\
 &\times\binom{q-6}{\frac{q-9}{2}}^2\binom{q-6}{\frac{q-7}{2}}^2\l(\binom{q-6}{\frac{q-11}{4}}+\phi(d)^{-1}\binom{q-6}{\frac{q-11}{4}}\r)d^{\frac{(q-1)(3q-21)}{8}} \pmod {\mathbb F_q^{\times^2}}.\\
\end{align*}
 If  $\phi(d)=-1 $ or $\binom{q-6}{\frac{q-11}{4}}\equiv 0 \pmod q$ , $\widetilde T_{(q-11)/2}(d,q)\equiv 0 \pmod q$. Otherwise,
\begin{align*}
    &\binom{q-6}{\frac{q-11}{4}}+\phi(d)^{-1}\binom{q-6}{\frac{q-11}{4}}\\
    &=  2\frac{(q-6)!}{(\frac{q-11}{4})!(\frac{3q-13}{4})!}\\
    &= 2\frac{(q-6)!\times\frac{3q-9}{4}\frac{3q-5}{4}\cdots \frac{3q+7}{4} }{((\frac{q-11}{4})!)^2\times\frac{q-7}{4}\frac{q-3}{4}\cdots \frac{3q+7}{4}}\\
    &\equiv \frac{-7(q-1)!}{(-1)^{\frac{q+9}{4}}(\frac{q-7}{4}\frac{q-3}{4}\cdots \frac{q-1}{2})^2}  \\
    &\equiv 7(-1)^{\frac{q+5}{4}}(q-1)! \pmod {\mathbb F_q^{\times^2}}. \\
\end{align*}
In this case, by Lemma 2.3 and $\phi(d)=1$ we obtain that
\begin{align*}
\widetilde T_{(q-11)/2}(d,q)&\equiv 7(-1)^{\frac{q+5}{4}}(q-1)!(-1)^{\frac{(q-3)(q-5)}{8}}d^{\frac{3q+3}{4}} \ or \ 0  \\
& \equiv 7(q-1)!d^{\frac{q+1}{4}} \ or \ 0 \pmod {\mathbb F_q^{\times^2}}. \\
\end{align*}
In view of the above, there is an element $y_q(d)\in \mathbb F_q$ such that
$$
\widetilde T_{(q-11)/2}(d,q)=7(q-1)!d^{\frac{q+1}{4}}y_q(d)^2.
$$ \qed

Applying this result, we can obtain Corollary 1.2 easily.

\noindent{\bf Proof of Corollary 1.2}.
 If $p\equiv 1 \pmod 4$ and $\l(\frac{d}{p}\r)=1$, from the proof above, we obtain
 $$\(\frac{T_{(p-11)/2}(d,p)}{p}\)= 1\ \t{or}\ 0.$$In other words, $\(\frac{T_{(p-11)/2}(d,p)}{p}\)\neq -1 .$

If $p\equiv 3 \pmod 4$ and $\l(\frac{d}{p}\r)=1$, we obtain that
$(\frac{T_{(p-11)/2}(d,p)}{p})= (\frac{-7}{p})$ or $0$.
 We get that if
 $$\(\frac{T_{(p-11)/2}(d,p)}{p}\)  = -1 ,$$
 then $p \equiv  1 \ \t{or} \ 2 \ \t{or} \ 4 \pmod 7$ immediately. \qed

\section{proof of Theorem 1.3}
\setcounter{lemma}{0}
\setcounter{theorem}{0}
\setcounter{equation}{0}
\setcounter{conjecture}{0}
\setcounter{remark}{0}
\setcounter{corollary}{0}
\noindent{\bf Proof } Let $a=(\frac {p-1}2)! .$ For $p\equiv 1 \pmod 4$ we have $a^2=((\frac {p-1}2)!)^2\equiv -1\pmod p$.
For each $k = 1, \cdots,(p-1)/2$, let
$\pi_a^*(k)$ be the unique $ r \in \{1, . . . ,(p-1)/2 \}$ with $ak$ congruent to $r$ or $-r$ modulo $p$. For the permutation $\pi_a^*$ on $\{1, \cdots ,(p-1)/2\}$,
Huang and Pan \cite{P06} showed that its sign is given by
$$ \mathrm{sgn}(\pi_a^*)=\l(\frac{a}{p}\r)^{(p+1)/2 }.$$
 If we view $D_p^{(m)}$ as a determinant over $\mathbb F_p$, then we have

$$D_p^{(m)}=\det \left[(i^2+(aj)^2)^{m+\frac{p-1}{2}}\right]_{1\leq i,j\leq (p-1)/2}=\mathrm{sgn}(\pi_a^*)\det \left[(a_i+a_j)^{m+\frac{p-1}{2}}\right]_{1\leq i,j\leq (p-1)/2}.$$
Thus
$$D_p^{(m)}\equiv \mathrm{sgn}(\pi_a^*)\det \left[(\frac {a_i}{a_j}+1)^{m+\frac{p-1}{2}}\right]_{1\leq i,j\leq (p-1)/2}\pmod {\mathbb F_p^{\times 2}}.$$
Fix positive odd integer $m$. For integer $k$ with $0 \leq k \leq (m-1)/2 $, we define $$F_m(k)=2^{m-2k}(m-k)(m-k-1)\ldots (k+1)+(2m-2k-1)(2m-2k-3)\ldots (2k+1).$$
Let
$$M=\max_{0 \leq k \leq (m-1)/2 }F_m(k).$$
We claim that $p\nmid D_p^{(m)}$ when $p>M$. In other words, if $p>M$, then we have $p \notin E(m).$ Hence, the elements of $E(m)$ is finite. Now, we prove the claim.

Noting that $p>M\geq F_m(0)=2^mm!+(2m-1)!!>2m$ and
$\l(\frac {a_i}{a_j}\r)^{\frac {p-1}2}=1, $

we can verify that
\begin{align*}
&\l(\frac {a_i}{a_j}+1\r)^{m+\frac{p-1}{2}}\\
=&\sum_{k=0}^{m+\frac{p-1}{2}}\binom{m+\frac{p-1}{2}}{k}\l(\frac {a_i}{a_j}\r)^k\\
=&\sum_{k=0}^{m}\l(\binom{m+\frac{p-1}{2}}{k}+\binom{m+\frac{p-1}{2}}{\frac{p-1}{2}+k}\r)\l(\frac {a_i}{a_j}\r)^k
+\sum_{k=m+1}^{\frac{p-3}{2}}\binom{m+\frac{p-1}{2}}{k}\l(\frac {a_i}{a_j}\r)^{k}.\\
\end{align*}
Hence
$$
\l(\frac {a_i}{a_j}+1 \r)^{m+\frac{p-1}{2}}=f\l(\frac {a_i}{a_j}\r)
$$
where
$$f(T)=\sum_{k=0}^{m}\l(\binom{m+\frac{p-1}{2}}{k}+\binom{m+\frac{p-1}{2}}{\frac{p-1}{2}+k}\r)T^k
+\sum_{k=m+1}^{\frac{p-3}{2}}\binom{m+\frac{p-1}{2}}{k}T^{k}.
$$
Let $C_f$ denote the product of coefficients of $f(T),$
then
$$
C_f=\prod_{k=0}^m\l(\binom{\frac{p-1}{2}+m}{k}+\binom{\frac{p-1}{2}+m}{m-k}\r)\prod_{k=m+1}^{\frac {p-3}2}\binom{\frac{p-1}{2}+m}{k}
$$
Applying Lemma 2.1 to $P(T)=f(T)$, we have
\begin{align*}
D_p^{(m)}&\equiv \mathrm{sgn}(\pi_a^*)\mathrm{sgn}(\mathrm{inv_p})\prod_{k=0}^m\l(\binom{\frac{p-1}{2}+m}{k}+\binom{\frac{p-1}{2}+m}{m-k}\r)\prod_{k=m+1}^{\frac {p-3}2}\binom{\frac{p-1}{2}+m}{k}\prod_{1\leq i<j\leq n}(a_j-a_i)^2\\
&\equiv \mathrm{sgn}(\pi_a^*)\mathrm{sgn}(\mathrm{inv_p})\prod_{k=0}^m\l(\binom{\frac{p-1}{2}+m}{k}+\binom{\frac{p-1}{2}+m}{m-k}\r)\prod_{k=m+1}^{\frac {p-3}2}\binom{\frac{p-1}{2}+m}{k}\pmod{\mathbb F_p^{\times^2}}.
\end{align*}
Since $p\nmid \prod_{m+1}^{\frac {p-3}2}\binom{\frac{p-1}{2}+m}{k} $ for $p>2m$, we can deduce that if $p\mid D_p^{(m)}, $ then
$$p \mid  \prod_{k=0}^m\l(\binom{\frac{p-1}{2}+m}{k}+\binom{\frac{p-1}{2}+m}{m-k}\r).  $$
For $0 \leq k \leq m $,
\begin{align*}
   &\binom{\frac{p-1}{2}+m}{k}+\binom{\frac{p-1}{2}+m}{m-k} =\frac{(\frac{p-1}{2}+m)!\l((m-k)!(\frac{p-1}{2}+k)!+k!(\frac{p-1}{2}+m-k)!\r)}{(m-k)!(\frac{p-1}{2}+k)!k!(\frac{p-1}{2}+m-k)!}.
\end{align*}
Hence, if $p \mid \binom{\frac{p-1}{2}+m}{k}+\binom{\frac{p-1}{2}+m}{m-k}$, we have
$$(m-k)!(\frac{p-1}{2}+k)!+k!(\frac{p-1}{2}+m-k)!\equiv 0 \pmod p .$$
Without loss of generality, we assume $k\leq (m-1)/2$ since $m$ is an odd positive integer. Then we can deduce that
\begin{align*}
p \mid (m-k)(m-k-1)\ldots (k+1)+(\frac {p-1}2+m-k)(\frac {p-1}2+m-k-1)\ldots (\frac {p-1}2+k+1).
\end{align*}
It means that
$$2^{m-2k}(m-k)(m-k-1)\ldots (k+1)+(2m-2k-1)(2m-2k-3)\ldots (2k+1)=F_m(k) \equiv 0 \pmod p, $$
which is in contradiction to $p>M.$ Therefore, $p \nmid \binom{\frac{p-1}{2}+m}{k}+\binom{\frac{p-1}{2}+m}{m-k}. $ Hence the claim is true.

We can also compute $E(m)$ for some specific $m$. For example, let $m=5$. Firstly, we check all the prime number $p$ with $p\leq 2m+3=13.$ We find that there is no prime $p\in E(5).$ For prime number $p>13,$ we get that if $p\mid D_p^{(5)}$, then $p \mid \binom{\frac{p-1}{2}+5}{k}+\binom{\frac{p-1}{2}+5}{5-k} $ for some $k$ from the illustration above. Moreover, we have $p\equiv 1 \pmod 4$
 and $p\mid F_5(k)$  for some $k$. We calculate $F_5(0),F_5(1)$ and $F_5(2)$ and only find that $29\mid F_5(0)=4785.$ Therefore, $E(5)=\{29\}.$

In the same way, we obtain that
$$ E(7)=\{13, 53\},\ E(9)=\{13, 17, 29\} $$
$$E(11)=\{17, 29\} \ and \ E(13)=\{17, 109, 401\}. $$ \qed

\section{proofs of Theorem 1.4 and 1.5}
\setcounter{lemma}{0}
\setcounter{theorem}{0}
\setcounter{equation}{0}
\setcounter{conjecture}{0}
\setcounter{remark}{0}
\setcounter{corollary}{0}
\noindent{\bf Proof of Theorem 1.4. }
By Cayley's Theroem we have $D_p^{(1)}=\mathrm{pf}(D_p^{(1)})^2$. For $p\equiv 1 \pmod 4$ we have $((\frac {p-1}2)!)^2\equiv -1\pmod p$. Let $a=(\frac {p-1}2)! $ and define $\pi_a^*(k)$ the same as before. By Lemma \cite[Lemma 2.3]{S19}, we have
$$\mathrm{sgn} (\pi_a^*)=\l(\frac{a}{p}\r)^{(p+1)/2 }=\l(\frac{2}{p}\r)^{(p+1)/2}=(-1)^{(p-1)/4}.$$
If we view $D_p^{(1)}$ as a determinant over $\mathbb F_p$,  we have
\begin{align*}
   D_p^{(1)}&=\det \left[(i^2+(aj)^2)^{\frac{p+1}{2}}\right]_{1\leq i,j\leq (p-1)/2}\\
   &=\mathrm{sgn}(\pi_a^*)\det \left[(a_i+a_j)^{\frac{p+1}{2}}\right]_{1\leq i,j\leq (p-1)/2}\\
 &=(-1)^{(p-1)/4}\prod_{i=1}^{(p-1)/2}a_i\det \left[\l(\frac {a_i}{a_j}+1\r)^{\frac{p+1}{2}}\right]_{1\leq i,j\leq (p-1)/2}\\
 &=(-1)^{(p+3)/4}\det \left[\l(\frac {a_i}{a_j}+1\r)^{\frac{p+1}{2}}\right]_{1\leq i,j\leq (p-1)/2}.\\
\end{align*}
Noting that
$\l(\frac {a_i}{a_j}\r)^{\frac {p-1}2}=1, $
 we can verify that
 \begin{align*}
   \l(\frac {a_i}{a_j}+1\r)^{\frac{p+1}{2}}= \sum_{k=2}^{\frac{p-3}{2}}\binom{\frac{p+1}{2}}{k}\l(\frac{a_i}{a_j}\r)^k+\frac{p+3}{2}\l(\frac{a_i}{a_j}\r)+\frac{p+3}{2}.\\
 \end{align*}

Hence for
$$f(T)=\sum_{k=2}^{\frac{p-3}{2}}\binom{\frac{p+1}{2}}{k}T^k+\frac{p+3}{2}T+\frac{p+3}{2}
$$
we have $$
\l(\frac {a_i}{a_j}+1 \r)^{\frac{p+1}{2}}=f\l(\frac {a_i}{a_j}\r).
$$
By Lemma 2.1 and Lemma 2.3, we obtain
 \begin{align*}
    D_p^{(1)}&=(-1)^{(p+3)/4}(-1)^{\frac{(p-3)(p-5)}{8}}\l(\frac{p+3}{2}\r)^2 \prod_{k=2}^{\frac{p-3}{2}}\binom{\frac{p+1}{2}}{k} \prod_{1\leq i<j\leq \frac{p-1}{2}}(i^2-j^2)^2\\
    &=\l(\frac{p+3}{2}\r)^2 \prod_{k=2}^{\frac{p-1}{4}}\binom{\frac{p+1}{2}}{k}^2 \prod_{1\leq i<j\leq \frac{p-1}{2}}(i^2-j^2)^2.
 \end{align*}
 Let
 $$t=\frac{p+3}{2} \prod_{k=2}^{\frac{p-1}{4}}\binom{\frac{p+1}{2}}{k} \prod_{1\leq i<j\leq \frac{p-1}{2}}(j^2-i^2).$$
 We have
 $$ D_p^{(1)}=\mathrm{pf}(D_p^{(1)})^2\equiv t^2
\pmod p. $$
Therefore $\sqrt{ D_p^{(1)}}=\pm t \pmod p.$ For $p\equiv 1 \pmod 4$ and $(\frac{-1}{p})=1,$ we obtain that
$$ \l(\frac{\sqrt{D_p^{(1)}}}{p}\r)=\l(\frac{t}{p}\r).$$
Referring to \cite{S20}, we have
$$\prod_{1\leq i<j\leq \frac{p-1}{2}}(j^2-i^2)=-(\frac {p-1}2)! \pmod p.
$$
Observe that
$$\prod_{k=2}^{\frac{p-1}{4}}\binom{\frac{p+1}{2}}{k}=\l((\frac{p+1}{2})!\r)^{\frac{p-5}{4}}\prod_{k=2}^{\frac{p-3}{2}}k!.
$$
By observation, we have
$$ \prod_{k=2}^{\frac{p-3}{2}}k!\equiv \frac{p-3}{2}(\frac{p-5}{2})^2(\frac{p-7}{2})^3\cdots 3^{\frac{p-3}{2}}2^{\frac{p-5}{2}} \equiv \frac{p-3}{2}\frac{p-7}{2}\cdots 3\equiv(\frac{p-3}{2})!!\pmod {\mathbb F_p^{\times2} }.$$
We also have
\begin{align*}
  (\frac{p-3}{2})!!&\equiv(\frac{p-3}{2})!!\l((\frac{p-1}{2})!!\r)^2  \\
  &\equiv (\frac{p-1}{2})! (\frac{p-1}{4})!2^{(p-1)/4} \pmod {\mathbb F_p^{\times2} }.
\end{align*}
We conclude that
 \begin{align*}
     t&\equiv -\frac{3}{2}\l((\frac{p+1}{2})!\r)^{\frac{p-5}{4}}(\frac{p-1}{2})! (\frac{p-1}{4})!2^{(p-1)/4}(\frac {p-1}2)!\\
     &\equiv -3\l((\frac{p+1}{2})!\r)^{\frac{p-5}{4}}(\frac{p-1}{4})!2^{(p-5)/4} \pmod {\mathbb F_p^{\times2} }.\\
 \end{align*}
For $p\equiv 1 \pmod 4$, by Lemma \cite[Lemma 2.3]{S19} we have
$$\l(\frac{(\frac{p+1}{2})!}{p}\r)=\l(\frac{2^{-1}}{p}\r)\l(\frac{(\frac{p-1}{2})!}{p}\r)=\l(\frac{2^{-1}}{p}\r)\l(\frac{2}{p}\r)=1 .$$
Notice that $\l(\frac{2^{(p-5)/4}}{p}\r)=1$ and $(\frac{-1}{p})=1$, we obtain
$$\l(\frac{\sqrt{D_p^{(1)}}}{p}\r)=\l(\frac{t}{p}\r)=\l(\frac{3}{p}\r)\l(\frac{(\frac{p-1}{4})!}{p}\r)=(-1)^{|\{0<k<\frac p4: (\frac kp)=-1\}|}\l(\frac p3\r).\qed
$$

\noindent{\bf Proof of Theorem 1.5. } For $p=5$, we check that Theorem 1.5 is right. Now, we consider $p\geq 13$ with $p\equiv 1 \pmod 4.$
By Cayley's Theorem we have $D_p^{(3)}=\mathrm{pf}(D_p^{(3)})^2$.

For $p\equiv 1 \pmod 4$ we have $((\frac {p-1}2)!)^2\equiv -1\pmod p$. Let $a=(\frac {p-1}2)! $ and define $\pi_a^*(k)$ the same as before. By Lemma \cite[Lemma 2.3]{S19}, we have
$$ \mathrm{sgn}(\pi_a^*)=\l(\frac{a}{p}\r)^{(p+1)/2 }=\l(\frac{2}{p}\r)^{(p+1)/2}=(-1)^{(p-1)/4}.$$
If we view $D_p^{(3)}$ as a determinant over $\mathbb F_p$, then we have
\begin{align*}
   D_p^{(3)}&=\det \left[(i^2+(aj)^2)^{\frac{p+5}{2}}\right]_{1\leq i,j\leq (p-1)/2}\\
   &=\mathrm{sgn}(\pi_a^*)\det \left[(a_i+a_j)^{\frac{p+5
   }{2}}\right]_{1\leq i,j\leq (p-1)/2}\\
 &=(-1)^{(p-1)/4}\prod_{i=1}^{(p-1)/2}a_i\det \left[\l(\frac {a_i}{a_j}+1\r)^{\frac{p+5}{2}}\right]_{1\leq i,j\leq (p-1)/2}\\
 &=(-1)^{(p+3)/4}\det \left[\l(\frac {a_i}{a_j}+1\r)^{\frac{p+5}{2}}\right]_{1\leq i,j\leq (p-1)/2}.\\
\end{align*}
Noting that
$\l(\frac {a_i}{a_j}\r)^{\frac {p-1}2}=1, $
 we can verify that
 \begin{align*}
   \l(\frac {a_i}{a_j}+1\r)^{\frac{p+5}{2}}=& \sum_{k=4}^{\frac{p-3}{2}}\binom{\frac{p+5}{2}}{k}\l(\frac {a_i}{a_j}\r)^k+1+\binom{\frac{p+5}{2}}{1}\l(\frac {a_i}{a_j}\r)+\binom{\frac{p+5}{2}}{2}\l(\frac {a_i}{a_j}\r)^{2}+\binom{\frac{p+5}{2}}{3}\l(\frac {a_i}{a_j}\r)^{3}\\
   &+\binom{\frac{p+5}{2}}{\frac{p-1}{2}}\l(\frac {a_i}{a_j}\r)^{\frac{p-1}{2}}+\binom{\frac{p+5}{2}}{\frac{p+1}{2}}\l(\frac {a_i}{a_j}\r)^{\frac{p+1}{2}}+\binom{\frac{p+5}{2}}{\frac{p+3}{2}}\l(\frac {a_i}{a_j}\r)^{\frac{p+3}{2}}+\binom{\frac{p+5}{2}}{\frac{p+5}{2}}\l(\frac {a_i}{a_j}\r)^{\frac{p+5}{2}} \\
   =&\sum_{k=4}^{\frac{p-3}{2}}\binom{\frac{p+5}{2}}{k}\l(\frac {a_i}{a_j}\r)^k+1+\binom{\frac{p+5}{2}}{\frac{p-1}{2}}+\l(\binom{\frac{p+5}{2}}{1}+\binom{\frac{p+5}{2}}{\frac{p+1}{2}}\r)\l(\frac {a_i}{a_j}\r) \\
   &+\l(\binom{\frac{p+5}{2}}{2}+\binom{\frac{p+5}{2}}{\frac{p+3}{2}}\r)\l(\frac {a_i}{a_j}\r)^2+\l(\binom{\frac{p+5}{2}}{3}+\binom{\frac{p+5}{2}}{\frac{p+5}{2}}\r)\l(\frac {a_i}{a_j}\r)^3\\
   =&\sum_{k=4}^{\frac{p-3}{2}}\binom{\frac{p+5}{2}}{k}\l(\frac {a_i}{a_j}\r)^k+\frac{21}{16}+\frac{35}{8}\l(\frac {a_i}{a_j}\r)+\frac{35}{8}\l(\frac {a_i}{a_j}\r)^2+\frac{21}{16}\l(\frac {a_i}{a_j}\r)^3.
 \end{align*}
Hence
$$
\l(\frac {a_i}{a_j}+1 \r)^{\frac{p+1}{2}}=f\l(\frac {a_i}{a_j}\r)
$$
where
$$f(T)=\sum_{k=4}^{\frac{p-3}{2}}\binom{\frac{p+5}{2}}{k}T^k+\frac{21}{16}+\frac{35}{8}T+\frac{35}{8}T^2+\frac{21}{16}T^3.
$$
By Lemma 2.1 and Lemma 2.3, we obtain
 \begin{align*}
    D_p^{(3)}&=(-1)^{(p+3)/4}(-1)^{\frac{(p-3)(p-5)}{8}}\l(\frac{21}{16}\r)^2\l(\frac{35}{8}\r)^2 \prod_{k=4}^{\frac{p-3}{2}}\binom{\frac{p+5}{2}}{k} \prod_{1\leq i<j\leq \frac{p-1}{2}}(i^2-j^2)^2\\
    &=\l(\frac{21}{16}\r)^2\l(\frac{35}{8}\r)^2 \prod_{k=4}^{\frac{p+3}{4}}\binom{\frac{p+5}{2}}{k}^2 \prod_{1\leq i<j\leq \frac{p-1}{2}}(i^2-j^2)^2.
 \end{align*}
 Let
 $$s=\frac{21}{16}\times\frac{35}{8} \prod_{k=4}^{\frac{p+3}{4}}\binom{\frac{p+5}{2}}{k} \prod_{1\leq i<j\leq \frac{p-1}{2}}(j^2-i^2).$$
 We have
 $$ D_p^{(3)}=\mathrm{pf}(D_p^{(3)})^2\equiv s^2
\pmod p. $$
Therefore $\sqrt{ D_p^{(3)}}=\pm s \pmod p.$ For $p\equiv 1 \pmod 4$ and $(\frac{-1}{p})=1,$ we obtain that
$$\l(\frac{\sqrt{D_p^{(3)}}}{p}\r)=\l(\frac{s}{p}\r).$$
Referring to \cite{S20}, we have
$$\prod_{1\leq i<j\leq \frac{p-1}{2}}(j^2-i^2)=-(\frac {p-1}2)! \pmod p.
$$
Observe that
$$\prod_{k=4}^{\frac{p+3}{4}}\binom{\frac{p+5}{2}}{k}=\l((\frac{p+5}{2})!\r)^{\frac{p-9}{4}}\prod_{k=4}^{\frac{p-3}{2}}k!.
$$
We notice that $\frac{p-3}{2}$ appears once,  $\frac{p-5}{2}$ appears twice ...  and $4$ appears $\frac{p-9}{2}$ times in
$\prod_{k=4}^{\frac{p-3}{2}}k!$.
Hence,
$$ \prod_{k=2}^{\frac{p-3}{2}}k!\equiv \frac{p-3}{2}\frac{p-7}{2}\cdots 5\equiv\frac{1}{3}(\frac{p-3}{2})!! \pmod {\mathbb F_p^{\times2} }.$$
We also have
\begin{align*}
  (\frac{p-3}{2})!!&\equiv (\frac{p-1}{2})! (\frac{p-1}{4})!2^{(p-1)/4} \pmod {\mathbb F_p^{\times2} }.
\end{align*}
We conclude that
 \begin{align*}
     s&\equiv -\frac{21}{16}\times\frac{35}{8}\l((\frac{p+5}{2})!\r)^{\frac{p-9}{4}}\frac{1}{3}(\frac{p-1}{2})! (\frac{p-1}{4})!2^{(p-1)/4}(\frac {p-1}2)!\\
     &\equiv -5\l((\frac{p+5}{2})!\r)^{\frac{p-5}{4}}(\frac{p-1}{4})!2^{(p-5)/4} \pmod {\mathbb F_p^{\times2} }.\\
 \end{align*}
For $p\equiv 1 \pmod 4$, by Lemma \cite[Lemma 2.3]{S19} we have
$$\l(\frac{(\frac{p+5}{2})!}{p}\r)=\l(\frac{{\frac{15}{4}}}{p}\r)\l(\frac{(\frac{p+1}{2})!}{p}\r)
=\l(\frac{3}{p}\r)\l(\frac{5}{p}\r) .$$
Notice that $\l(\frac{2^{(p-5)/4}}{p}\r)=1$ and $(\frac{-1}{p})=1$, we obtain
\begin{align*}
    \l(\frac{\sqrt{D_p^{(3)}}}{p}\r)&=\l(\frac{s}{p}\r) \\
 &=\l(\frac{3}{p}\r)^{(p-9)/4}\l(\frac{5}{p}\r)^{(p-5)/4}\l(\frac{(\frac{p-1}{4})!}{p}\r)\\
 &=(-1)^{|\{0<k<\frac p4: (\frac kp)=-1\}|}\l(\frac {p}{4+(-1)^{(p-1)/4}}\r).\\
\end{align*} \qed

\Ack The authors would like to thank Zhi-Wei Sun for  providing us these conjectures as research project in January 2024.

\end{document}